\documentclass[11pt]{amsart}
\usepackage{amscd}
\renewcommand{\baselinestretch}{1.1}

\newtheorem{theorem}{Theorem}[section]
\newtheorem{lemma}[theorem]{Lemma}
\newtheorem{corollary}[theorem]{Corollary}
\newtheorem{proposition}[theorem]{Proposition}
\newtheorem{hypothesis}{Hypothesis}
\theoremstyle{definition}
\newtheorem{definition}[theorem]{Definition}

\theoremstyle{remark}

\numberwithin{equation}{section}

\DeclareMathOperator{\End}{End}
\newcommand{\Nab}[2]{\nabla^{#1}_{#2}}
\newcommand{\Cinf}{C^{\infty}}

\newcommand{\SF}{SF}

\def\span{\hbox{span}}
\def\proj{\hbox{proj}}
\def\proof{\noindent{\sl Proof.\ \ }}

\def\text#1{{\mbox{#1}}}
\def\norm{\|}

\def\al{\alpha}

\def\la{\lambda}

\def\si{\sigma}
\def\Ga{\Gamma}
\def\La{\Lambda}

\def\Si{\Sigma}
\def\cP{{\mathcal P}}

\def\wtilde{\widetilde}

\def\lto{\longrightarrow}

\renewcommand{\span}{\operatorname{span}}

\def\sqr#1#2{{\vcenter{\hrule height.#2pt 
\hbox{\vrule width.#2pt height #1pt \kern#1pt \vrule width.#2pt}
\hrule height.#2pt}}}
\def\square{\mathchoice\sqr{5.5}4\sqr{5.0}4\sqr{4.8}3\sqr{4.8}3}
\def\qed{\hskip4pt plus1fill\ $\square$\par\medbreak\vskip5ex}

\begin{document}
\begin{abstract}
We derive a decomposition formula for the spectral flow of a 1-parameter family
of self-adjoint Dirac operators on an odd-dimensional
manifold $M$ split along a hypersurface $\Sigma$ ($M=X\cup_{\Sigma} Y$).  No
transversality or stretching hypotheses are assumed and the boundary conditions
can be chosen arbitrarily.   The formula takes the form $SF(D)= SF(D_{|X}, B_X)
+ SF(D_{|Y},B_Y) + \mu(B_Y,B_X) + S$ where $B_X$ and $B_Y$ are
boundary conditions, $\mu$ denotes the Maslov index, and
$S$ is a sum of explicitly defined Maslov indices coming from stretching and
rotating boundary conditions. The derivation is a simple consequence of
Nicolaescu's theorems and elementary properties of the Maslov index.  We show
how to use the formula and derive many of the splitting theorems in the
literature as simple consequences.
\end{abstract}

\title[A general splitting formula . . .] {A general splitting formula for the
spectral flow}

\author{M. Daniel and P. Kirk, with an appendix by K.P. Wojciechowski}
\date{\today}
\maketitle
\section{Introduction}

Several articles have been written containing formulas expressing the spectral
flow of a path of self-adjoint Dirac operators on a closed, split manifold $M$
($M=X\cup_{\Sigma} Y$)  in terms of quantities determined by each piece in the
decomposition and ``interaction'' terms. For example see \cite{bunke, CLM2,
Yoshida, Nico, AMD}. The article of Nicolaescu  \cite{Nico} is perhaps  the most
elegant  and conceptually   appealing.  Additionally, a large number  of articles
consider the closely related but more delicate problem of splitting theorems for
the Atiyah-Patodi-Singer invariant. The bibliography to Bunke's article
\cite{bunke} contains a long list of citations.  Most of these articles, with the
exception of Nicolaescu's, use delicate analytical methods and estimates such as
heat kernel methods, and the results apply only after one has stretched the collar
neighborhood of the separating hypersurface.   Nicolaescu instead treats the
problem largely from the point of view of linear algebra in a symplectic Hilbert
space, and his main result is appealing in the simplicity of its statement:  the
spectral flow of the path equals the Maslov index
$\mu(\Lambda_X,\Lambda_Y)$.  Here $\Lambda_X$ and $\Lambda_Y$ denote the  paths
of   Cauchy data spaces  consisting of the restrictions of nullspace elements of
the operators on $X$ and
$Y$   to their common boundary $\Sigma$.  Unfortunately, Nicolaescu's
formulation does not lend itself easily to computation. What is needed  is a
splitting formula that isolates the contribution from each of  the two
pieces of the decomposition to the spectral flow.  This is especially important
when studying spectral flow in the context of cut-and-paste constructions.

 In this article we prove a general splitting theorem and show how it can be used
to derive most of the various splitting theorems in the literature. The proof of
our result  is quite simple, and uses only elementary properties of the Maslov
index in addition to three results of
 Nicolaescu: the theorem described in the preceding paragraph, a version from his
subsequent article for manifolds with boundary \cite{Nico2}, and the calculation
of the adiabatic limit of the  Cauchy data space from \cite{Nico}.

Our main result, Theorem
\ref{splitthm} states:
\vskip3ex

\noindent{\bf Theorem.} {\it  Let $D(t)$ be a continuous path of self-adjoint Dirac
operators on a smooth, closed, oriented, odd-dimensional, Riemannian manifold
$M$.  Suppose that $M$ can be split along a hypersurface $\Si$ ($M=X\cup_{\Si}Y$)
and that each
$D(t)$ is cylindrical  and neck-compatible with respect to this splitting.  Let
$B_X(t)$ and
$B_Y(t)$ be paths of self-adjoint elliptic boundary conditions for the restriction
of
$D(t)$ to $X$ and $Y$ respectively.   

Then 
\[  SF(D)= SF(D_{|X},B_X) + SF(D_{|Y},B_Y) +\mu(B_Y(1-t),B_X(1-t))
+\sum_{i=1,2,4,5,7,8,10,11} \mu(L_i,M_i)\] }

The terms appearing in the sum are certain Maslov indices of explicitly defined
paths of Lagrangians. 

Notice that this formula, in contrast to the theorems cited above, holds without
any preliminary stretching assumptions, nor any prescription on what the boundary
conditions should be.

Perhaps the method itself is more important than the actual formula, in the sense
that in any given application it is probably easier to adapt the method we
introduce here to the specific situation than to make the problem fit our formula.
(This is the case in the article
\cite{BHKK} on the $SU(3)$ Casson invariant.)  For that reason we include a
lengthy ``user's guide'' (Section \ref{UsersGuide}) which indicates how various
additional hypotheses can be used to force some of the terms $\mu(L_i,M_i)$ to
vanish.  We also show how to easily derive many of the different versions of the
splitting theorems cited above. In particular, we derive the splitting theorem 
of Bunke,  give a generalization of this theorem and the splitting theorem of
Yoshida and Nicolaescu, and indicate the relation between our formula and the
formula of
\cite{CLM2}.

Our results are stated and proven for Dirac operators on odd dimensional
manifolds since these include most of the geometrically important classes of
self-adjoint elliptic operators such as the odd signature operator and the spin
Dirac operator.   
\vskip3ex

We finish this introduction with a brief example of the method  for those readers
who are familiar with  this subject. Other readers    can return to the following
paragraphs after  finishing Section
\ref{four}.

Suppose that $D(t):\Ga(E)\to\Ga(E)$, $t\in [0,1]$  is a path of self-adjoint Dirac
operators on a manifold $M$ decomposed along a hypersurface $M=X\cup_\Si Y$.  Let
$\Lambda_X(t)$ and $\La_Y(t)$ be the Cauchy data spaces associated to the
restrictions of $D(t)$ to $X$ and $Y$ respectively. These are Lagrangian subspaces
of the symplectic Hilbert space $L^2(E|_{\Si})$.    Assume further that each
$D(t)$ is cylindrical ($D(t)=J(\partial/\partial s + S(t))$ on a collar
neighborhood of $\Si$) and neck-compatible (for each $t$,
$S(t)$ is self-adjoint). Furthermore, suppose that the kernels of the tangential
operators
$S(t)$ are trivial for all $t$ and denote by  $P^\pm(t)$ the positive/negative
eigenspace of
$S(t)$.  

Finally, suppose that $\La_X(0)=P^-(0)$, $\La_X(1)=P^-(1)$, $\La_Y(0)=P^+(0)$, and 
$\La_Y(1)=P^+(1)$.  These four equalities  rarely hold except in artificial
examples, but   Nicolaescu's adiabatic limit theorem says these conditions are
asymptotically true;  compensating for this leads to the extra terms in our
formula.

The path $\La_X(t)$ is clearly homotopic rel endpoints to the composite of
the three paths 
$P^-(t)$, $P^-(1-t)$ and $\La_X(t)$. Similarly the path $\La_Y(t)$ is  
homotopic rel endpoints to the composite of the three paths $\La_Y(t)$,
$P^+(1-t)$, and
$P^+(t)$.  Because the Maslov index is invariant under rel endpoint homotopies and
additive with respect to compositions of paths, we conclude
\begin{eqnarray*} SF(D,M)&=&\mu(\La_X,\La_Y)\ \ \  \hbox{ (Nicolaescu's splitting
theorem) }\\
 &=& \mu(P^-,\La_Y)  +\mu(P^-(1-t), P^+(1-t)) + \mu(\La_X; P^+)\\ &=& SF(D_{|Y};
P^-)  + SF(D_{|X};P^+)\end{eqnarray*} The last step follows from the version of
Nicolaescu's theorem for manifolds with boundary and the fact that $P^+$ and $P^-$
are transverse.  The proof of our main result is no more difficult than this. The
extra terms come about by moving to the adiabatic limits at the endpoints and from
allowing general boundary conditions.

The authors thank H. Boden, D. Hoff,  K. P. Wojciechowski
for helpful discussions, and especially L. Nicolaescu who first gave us a proof
of Lemma \ref{stretchpath}.  The authors thank K. P. Wojciechowski
 for letting us use his proof of this lemma in the appendix.

\section{Dirac operators}\label{DiracOperators} There are many different
definitions of {\em Dirac operator} in the literature. For our purposes, we adopt
that of \cite{Nico}.  Briefly, a  Dirac operator is determined by a Clifford
module over a manifold along with a  compatible connection.  More precisely,
suppose we are given the following
\begin{enumerate}
\item An oriented Riemannian manifold $(M,g)$.
\item A {\em self-adjoint Clifford module} $E\to M$.  So $E$ is a vector  bundle
over $M$ with an action $c:C(M)\to \End(E)$.  Here $C(M)$ is the  bundle of
Clifford algebras over
$M$ generated by the cotangent bundle using  the metric.  The adjective
self-adjoint means that $c$ carries each element  of $T^*M$ to a {\em
skew-adjoint} endomorphism.  Together with the Clifford relations, this implies
that the each element of $T^*M$ acts orthogonally.  For convenience we assume the
vector bundle $E$ is a complex vector bundle.
\item A {\em Clifford compatible covariant derivative} $\Nab{E}{}$ on
$E$.  Thus
\[\Nab{E}{}:\Gamma(E)\to\Gamma(E\otimes T^*M)\] is a differential operator
satisfying the Leibnitz rule
\[\Nab{E}{}(fs) = df\otimes s + f\Nab{E}{}s\] for any $f\in\Cinf(M)$ and
$s\in\Gamma(E)$, and compatible with the  Clifford action in the sense that 
\[ [\Nab{E}{},c(a)] = c(\Nab{}{}a) \] where $a\in \Gamma(C(M))$ and $\Nab{}{}$ is
the Levi-Civita connection  (naturally extended from $TM$ to $C(M)$).
\end{enumerate} This data determines a Dirac operator as the composition
\[
\begin{CD}
\Gamma(E) @>\Nab{E}{}>> \Gamma(E\otimes T^*M) @>CC>> \Gamma(E)
\end{CD}
\] where $CC$ denotes contraction with respect to the Clifford action (denoted by
$c$ above).  Incidentally, this definition agrees with that of a {\em Dirac
operator on a Dirac bundle} as defined in
\cite{Spin}.  In this article, we consider only self-adjoint Dirac operators over
odd-dimensional manifolds.

We are particularly interested in Dirac operators over split manifolds.  A
manifold $M$ is {\em split along a hypersurface
$\Sigma$ } if it can be expressed as the union of two manifolds with boundary,
$X$ and $Y$, such that $\partial X = -\partial Y = \Sigma = X\cap Y$.  In this
case we also require the existence of a neighborhood $U$ of $\Sigma$ in $M$ 
that is isometric to
$\Sigma\times (-1,1)$.  Over this neighborhood, all relevant structures ({\em
e.g.} the Clifford bundle $E$)  should decompose similarly.

Thus we are led to consider Dirac operators on manifolds with boundary, and in
this context we impose two further restrictions on such operators.  Such a Dirac
operator must be {\em cylindrical}, meaning that in a neighborhood of the boundary
(of the form $\Sigma\times (-1,0]$ or $\Sigma\times [0,1)$ as described above)
$D$ can be written as 
\begin{equation}\label{cylind}
 D = c(du)(\partial/\partial u + S) 
\end{equation} where $u$ is the second factor in $\Sigma\times (-1,0]$ (or
$\Sigma\times [0,-1)$), chosen so that $||du||=1$,  and $S$ is a Dirac operator
on
$E|_{\partial M}$, referred to as the {\em tangential operator}.   Note that $S$
is assumed to be constant in that it does not depend on the coordinate $u$. 
Finally, we require that $D$ be {\em neck compatible}, meaning that the
tangential operator $S$ is self-adjoint. 

In what follows we consider only Dirac operators satisfying these conditions. 
Although these conditions may appear restrictive, most  important
geometrically defined self-adjoint operators are of this type {\em e.g.} the spin
Dirac and  odd signature operators.

The Clifford relation ($v\otimes w + w\otimes v = -2\langle v, w\rangle$) implies 
that the algebraic operator $c(du):\Gamma(E_{|\partial M})\to
\Gamma(E_{|\partial M})$  is a fiberwise isometry satisfying $c(du)^2=-$Id, and 
so it induces a complex structure on $L^2(E_{|\partial M})$ which we rename
suggestively
 \begin{equation}\label{complx}
J:L^2(E_{|\partial M})\to L^2(E_{|\partial M}). 
\end{equation}

So $J^2=-$Id. Moreover, $SJ=-JS$ and so the spectrum of the elliptic self-adjoint
operator
$S:L^2(E_{|\partial M})\to L^2(E_{|\partial M})$ is symmetric, and its
$\lambda$ and $-\lambda$ eigenspaces are interchanged by $J$. 

Define a hermitian symplectic structure on $L^2(E_{|\partial M})$ by 
$$\omega(x,y)=\langle x, Jy\rangle$$ where $\langle\ , \ \rangle$ denotes the
$L^2$ inner product. 

\begin{definition} Two closed subspaces $L_1,L_2$ of a Hilbert space form a {\it
Fredholm pair} if
$L_1\cap L_2$  is   finite dimensional, and $L_1+L_2$ is closed with finite
codimension.    
\end{definition}

\begin{definition} 
\begin{enumerate}
\item A closed subspace $L\subset L^2(E_{|\partial M})$ is called {\it isotropic}
if $L$ and 
 $JL$ are orthogonal. Thus $\omega(l,m)=0$ for all $l,m\in  L$.
 \item A closed subspace $L\subset L^2(E_{|\partial M})$ is called {\it
Lagrangian} if
 $JL$ is the orthogonal complement of $L$. Thus $\omega(l,m)=0$ for all $l,m\in 
L$ and
$L+JL=L^2(E_{|\partial M})$.  
\end{enumerate}
\end{definition}

Since $SJ=-JS$, the Hilbert space $L^2(E_{|\partial M})$  has an orthogonal
decomposition into the orthogonal direct sum of the negative eigenspace, kernel,
and positive eigenspace of $S$
\begin{equation}\label{decomp1}  L^2(E_{|\partial M})= P^-(S)\oplus \ker S\oplus
P^+(S).
\end{equation} In this decomposition, $\ker S$ is   finite dimensional since $S$
is elliptic on the closed manifold $\partial M$. Moreover $J$ preserves $\ker S$
and so $\ker  S$ is a symplectic subspace.  The spaces 
$P^+(S)$ and $P^-(S)$  are  interchanged by $J$ since $JS=-SJ$, and so $P^+(S)$ 
and
$P^-(S)$ are  infinite dimensional  and isotropic.  

If $L\subset \ker S$ is a (finite dimensional) Lagrangian subspace (this is
defined just as before, substituting $\ker S$ for $L^2(E_{|\partial M})$), then 
the spaces $P^-(S)\oplus L$ and $L\oplus P^+(S)$ are easily seen to be Lagrangian
subspaces of
$L^2(E_{|\partial M})$.  An important case occurs when $\ker S=0$,  in which case
$P^\pm(S)$ are themselves Lagrangian subspaces.

It will be convenient to have a slightly more general decomposition of
$L^2(E_{|\partial M})$ than Equation \ref{decomp1}.  To this end, let
$\nu$  be any nonnegative real number and  define
\begin{equation}\label{hnu} H_\nu(S)=\span_{L^2}\{ \phi\ | \ S\phi=\lambda\phi
\hbox{ and } |\lambda|\leq
\nu\},
\end{equation}\label{pminmu}
\begin{equation} P^-_\nu(S)=\span_{L^2}\{ \phi\ | \ S\phi=\lambda\phi \hbox{ and
}  \lambda < -\nu\},
\end{equation}
 and 
 \begin{equation}\label{pplusnu} P^+_\nu(S)=\span_{L^2}\{ \phi\ | \
S\phi=\lambda\phi
\hbox{ and }  \lambda>
\nu\}.
\end{equation}  Then as before the $P^\pm_\nu(S)$ are infinite dimensional
isotropic subspaces  and
$H_\nu$ is a finite dimensional symplectic subspace.  Moreover the decomposition 
of Equation
\ref{decomp1} is a special case  ($\nu=0$)   of the   decomposition
\begin{equation}\label{decomp2} L^2(E_{|\partial M})= P^-_\nu(S)\oplus
H_\nu(S)\oplus P^+_\nu(S).
\end{equation}

It is proven in \cite{KK-illinois} that  if $S$ is taken to vary continuously
 over some parameter space $T$, i.e. if the map $t\mapsto S(t)-S(t_0)$ is a
continuous map from $T$ into the space of bounded operators (here $t_0$ is some
fixed base point in $T$) and $\nu(t)$ is a continuous non-negative function on 
$T$ so that   $S(t)$ has a spectral gap at $\nu(t)$ ({\em i.e.}
 $\nu(t)$ misses the spectrum of $S(t)$), then the decomposition (\ref{decomp2})
is continuous in $T$.  Continuity for subspaces will always be taken in the gap
topology
\cite{kato}.

\section{Cauchy data spaces}

For a given  Dirac operator $D$ on a manifold $X$ with non-empty boundary $\Si$, 
its {\em Cauchy data space}
$\Lambda_X(D)$ is a Lagrangian subspace of 
$L^2(E|_{\partial M})$ consisting roughly of boundary values of its kernel
elements.  We give a definition  suitable for our purposes,  referring to
\cite{Nico}  for a careful construction.

In \cite{BW} it is shown that in the present context there is a well defined,
bounded, injective   restriction map (see Proposition 2.2 of \cite{Nico})
\begin{equation}\label{restrictionmap} R:\ker\left(D:L^2_{{1\over2}}(E)\to
L^2_{-{1\over2}}(E) \right)
\lto  L^2(E_{|\Si}).
\end{equation} 
Here $L^2_s(E)$ means the Sobolev space of sections of $E$ with $s$ derivatives
in $L^2$, extended in the usual way to real $s$.

The image of $R$ is a closed, infinite dimensional
Lagrangian subspace of
$L^2(E_{|\Si})$.  It will be denoted by  
\begin{equation}\label{cauchy}
\La_X(D):= R\left(\ker\left(D:L^2_{{1\over2}}(E)\to L^2_{-{1\over2}}(E)
\right)\right)
\end{equation}
 and called the {\it Cauchy data space} of the operator $D$ on $X$.  Sometimes we
will  abbreviate $\La_X(D)$ to
$\La_X$ or even $\La$ when $D$ or $X$ are clear from context.  Thus the Cauchy
data space is space of  boundary values of solutions to $D\si=0$.   In
\cite{Nico} it is proven that if 
$D$ varies regularly (smooth is sufficient but not necessary) in the space of
Dirac operators with respect to some parameter space
$T$, then the Cauchy data spaces $\La_X(D(t))$ vary regularly (at least $C^1$) in
$t\in T$. Regularity for closed subspaces may be interpreted in terms of the norm
topology of the associated projections.  The resulting topology is equivalent to
the gap topology \cite{kato}.
\vskip3ex

An important property of the Cauchy data space of a Dirac operator $D$ of the form
$J(\partial /\partial u + S)$ on the collar $\Si\times[-1,0]$ of the boundary of
$X$ is that  the pair $(\La_X(D), P^+(S))$ forms a Fredholm pair of subspaces
\cite{Nico}.  Since $P^+_\nu(S)\subset P^+(S)$ has finite codimension, it follows 
that if
$B$ is any closed subspace of $L^2(E_{|\Si})$ which contains $P^+_\nu(S)$ for some
$\nu$ with finite codimension, then $(\La_X(D),B)$ form a Fredholm pair.

\vskip3ex

The proof of our main theorem will require stretching, which we now describe.  
Given a manifold $X$ with boundary $\Si$ and (open) collar $\Si\times (-1,0]$,
define, for $r\ge 0$, 
\begin{equation}\label{MR} X^r=X\cup_{\Si\times(-1,0]}\Si\times(-1,r].
\end{equation} Thus $X=X^0$.   Using Equation \ref{cylind}  to {\it define} $D$ on 
$\Si\times(-1,r]$ gives a natural extension of $D$ to $X^r$.  In this way one
obtains a 1-parameter family of Cauchy data spaces 
$\La_{X^r}(D)$.  The limit of $\La_{X^r}(D)$ as $r$ approaches infinity is
identified in Theorem 4.9 of \cite{Nico}.  We elaborate on this important and
interesting result.  

For notational convenience we write $\La^r_X$ for $\La_{X^r}(D)$ and $P^+_\nu$
for
$P^+_\nu(S)$.  Since $\La^0_X\cap P^+_0$ is finite dimensional, and since
$\cap_{\nu\to\infty}P^+_\nu=0$, there exists a number $\nu_0\ge 0$ so that
\begin{equation}\label{nonres}
\La^0_X\cap P^+_{\nu_0}=0.
\end{equation}  Following Nicolaescu, the set of all non-negative real numbers
satisfying  Equation \ref{nonres} is a non-empty, closed,  unbounded interval
called the  {\it nonresonance range of $D$}. The smallest such
$\nu_0$ is called the {\it nonresonance level of $D$}. Fix some  $\nu_0$ in the
nonresonance range of $D$. 

The symplectic reduction of $\La^0_X$ to $H_{\nu_0}$ is the Lagrangian subspace
\begin{equation}\label{eq1}
\wtilde{\La}_X(D)= \proj_{H_{\nu_0}}(\La^0_X\cap(H_{\nu_0}\oplus P^+_{\nu_0}))= {{
\La_X^0\cap(H_{\nu_0}\oplus P^+_{\nu_0}) }\over  { \La_X^0\cap P^+_{\nu_0}
}}\subset H_{\nu_0}.\end{equation}

The decomposition of Equation \ref{decomp2} is preserved by $S$ since this is an  
decomposition in terms of eigenspaces of $S$.  In particular $S$ preserves
$H_{\nu_0}$ and the restriction of $S$ to $H_{\nu_0}$  is self-adjoint with all
eigenvalues in
$[-\nu_0,\nu_0]$. Thus we can form the 1-parameter family of (finite dimensional)
operators 
\begin{equation}\label{eq2} e^{-rS}:H_{\nu_0}\to H_{\nu_0}.\end{equation}

It is not too hard to see that the limit
\begin{equation}\label{defofL} L_X(D):=\lim_{r\to \infty}
e^{-rS}\wtilde{\La}_X(D)
\end{equation} exists and is a Lagrangian subspace of $H_{\nu_0}$. 

We may now state Nicolaescu's adiabatic limit theorem \cite{Nico}.
\begin{theorem}\label{NicoAdiabaticLimitTheorem} As $r\to\infty$
\begin{equation}\label{adiabthm}
\La^r_X(D) \to  P^-_{\nu_0}\oplus L_X(D).
\end{equation}
\end{theorem} The limiting subspace is called the {\em adiabatic limit} of
$\La^r_X$.  Thus the adiabatic limit is determined, up to a finite dimensional
piece, by the tangential operator.

The identification of the adiabatic limit is an important ingredient in the proof
of our splitting formula, but we require a little bit more. We complement the
previous theorem with a lemma stating that the adiabatic deformation is in fact
regular.
\begin{lemma}\label{stretchpath} Let $r(t)={1\over{1-t}}$ for $t\in [0,1)$. The
path of Lagrangian subspaces
$$ t\mapsto \begin{cases}
\La^{r(t)}& t<1,\\ P^-_{\nu_0}\oplus  L_X(D) &t=1.
\end{cases}
$$  is continuous.
\end{lemma}
    The proof of Lemma \ref{stretchpath}  was provided to us by K.P.
Wojciechowski and can be found in the appendix.

 One warning is in order here. It is not true that the adiabatic limits of the
Cauchy data spaces vary continuously when 
$D$ is varying continuously over some parameter space, even if $\nu_0$ is larger
than the nonresonance level for every operator $D$  in this family.  The reason
for this is that the dynamics of $e^{-rS}$ acting on subspaces of $H_{\nu_0}$ is
quite sensitive the initial subspace.  See \cite{BHKK} for an explicit example of
an analytic path of Dirac operators
$D(t)$ for which the path of adiabatic limits of the  Cauchy data spaces
$P^-(t)\oplus L_X(D(t))$ is not continuous.  There are some special circumstances
when one can conclude that the adiabatic limits vary continuously, and in those
cases a splitting theorem can be proven easily. One such example is Theorem
\ref{yoshida} below.
\vskip3ex

We now let $M$ be a closed manifold split along a hypersurface $\Si$ into two
pieces
$X$ and $Y$ ($M=X\cup_{\Si}Y$). As above, we identify a closed neighborhood of
$\Si$ in $M$ as
$\Si\times(-1,1)$, with $\Si=\partial X=-\partial Y$. In the previous paragraphs
we have stated various facts about Dirac operators from the point of view of the
``$X$ side''. For convenience we state the analogous facts for the ``$Y$ side''.
The main thing to keep in mind here is that the complex structure $J$ on
$L^2(E|_{\Sigma})$ and the cylindrical decomposition (\ref{cylind}) in the collar
use the outward normal to
$X$ which is the inward normal to $Y$.  This generally has the effect of switching
signs.  Chasing down the repercussions one obtains the following facts (see
\cite{Nico}). 
\begin{enumerate}
 \item  $(P^-(S),\La_Y(D))$ is a Fredholm pair.
\item The limit as $r\to \infty$ of $\La_{Y^r}(D)$ is $L_Y(D)\oplus
P^+_{\nu_1}(S)$ where $\nu_1$ is in the nonresonance range of $D$ acting on $Y$
and 
$L_Y(D)$ is defined similarly to $L_X(D)$ but by taking $r\to -\infty$.
\item The pair $(\La_X(D),\La_Y(D))$ is a Fredholm pair.
\item The kernel of $D:\Ga(E)\to \Ga(E)$ is taken isomorphically to the
intersection 
$\La_X(D)\cap \La_Y(D)$ by restricting to $\Si$.
\end{enumerate}

\section{Spectral flow equals Maslov index}\label{four}

The theorems of Nicolaescu presented in this section establish the equality of two
a priori different invariants that can be associated to a path of Dirac
operators.  Accordingly, we begin with a description of these two invariants.

The spectrum of a Dirac operator $D$ on a closed manifold $M$ consists of discrete
eigenvalues of finite multiplicity.    The {\em spectral flow} of    a continuous 
path $D(t)$ ($t\in[0,1]$) of self-adjoint Dirac operators is (roughly) defined to
be the algebraic count (with multiplicity) of the number of eigenvalues crossing
through zero. While this definition is somewhat imprecise, it suffices for our
purposes, particularly because we never actually work with the spectral flow
directly.  Instead, we use Nicolaescu's theorems below to convert the spectral
flow to the Maslov index.  In any case, precise definitions of the spectral flow
can be found in \cite{KK, Nico, CLM2}.

An important technical point is appropriate here. One must set conventions so that
the  spectral flow is well defined on paths $D(t)$ for which
$D(0)$ and/or $D(1)$ have non-trivial kernel. One must decide
whether or not an eigenvalue that {\em starts} or {\em ends} at
$0$ counts as crossing {\em through} $0$.  It is important to be precise here,
because different conventions appear in the literature, and the particular choice
effects the properties of the invariant.  Among the several such conventions that
can be found in the literature, we will use the following.   Given a path $D(t),\
t\in[0,1]$ of Dirac operators, let $\epsilon>0$ be a number smaller than the
smallest positive eigenvalues of $D(0)$  and $D(1)$.  We {\em define} the spectral
flow of the path $D(t)$ to be the spectral flow of the path
$D(t)-\epsilon\hbox{Id}$: 
$$SF(D(t)):=SF(D(t)-\epsilon\hbox{Id}).$$  Effectively, we count the eigenvalues
that cross $\epsilon$ rather than those that cross $0$. Notice that this avoids
the issue of starting or ending at the crossing value, because by definition, no
eigenvalues start or end at $\epsilon$.

\vskip3ex

Given a continuous path  of Fredholm pairs of Lagrangians $(\La_1(t),\La_2(t))$ 
in a symplectic vector space the {\it Maslov index} $\mu(\La_1,\La_2)$ is the
integer defined to be the algebraic count of how many times $\La_1(t)$ passes
through $\La_2(t)$ along the path. The complex structure $J$ is used to specify
the signs in this algebraic count. In particular, the normalization is chosen so
that $\mu(\La_1,e^{tJ}\La_2), \ t\in [-\epsilon,\epsilon]$ equals
$\dim(\La_1,\La_2)$  when $\La_1$ and $\La-2$ are constant paths.

 See
\cite{Nico,AMD,CLM4} for the precise definition.  Note, the condition that the
Lagrangians be Fredholm is vacuous in the finite dimensional case, but critical in
our context ($L^2(E|_{\Sigma})$).  The Fredholm property is typically easily
verified for any pair of paths we consider by appealing to facts about Cauchy data
spaces and related Lagrangians as discussed in Section \ref{DiracOperators}.    

As with the spectral flow, a convention must be chosen to define the Maslov index
for paths of pairs that are not transverse at the endpoints. Again, it is
important to be explicit here because there are a number of possibilities. We
use a convention defined in terms of  the complex structure $J$ as explained in
\cite{CLM4}.  Choose a small positive $\epsilon$ so that 
\begin{enumerate}
 \item $(e^{s J}L_1(t), L_2(t))$ form a Fredholm pair for each $t$ and each $0\leq
s\leq\epsilon$. This is possible since Fredholm pairs form an open subspace of the
space of closed pairs \cite{kato}. 
\item $e^{sJ} L_1(0)$ is transverse to $L_2(0)$ and $e^{sJ} L_1(1)$ is transverse
to
$L_2(1)$ for all $0<s\leq \epsilon$.  The proof that such an $\epsilon$ exists can
be found in  \cite{CLM4}.
\end{enumerate}

Thus the path of pairs $(e^{J\epsilon}L_1(t),L_2(t))$ forms a path of Fredholm
pairs which are transverse at the endpoints. One then {\em defines} the Maslov
index of $L_1$ and $L_2$ by taking
\begin{equation}\label{nontrans}\mu(L_1,L_2):=\mu(e^{J\epsilon}L_1,L_2).
\end{equation}

We will use the following two elementary properties of the Maslov index.
\begin{enumerate}
\item {\em Path Additivity}:  Let $L_1$, $L_2$, $K_1$, and $K_2$ be paths of
Lagrangians such that $L_i(1) = K_i(0)$ for $i=1,2$ and let
$M_i$ be the path obtained by concatenating $L_i$ and $K_i$ (we write
$M_i=L_i*K_i$).  Then
\[\mu(M_1,M_2) = \mu(L_1,L_2) + \mu(K_1,K_2).\] 
\item{\em Homotopy Invariance}: Let $L_1$, $L_2$, $K_1$, and $K_2$ be paths of
Lagrangians such that $L_i$ is homotopic rel endpoints to
$K_i$.  Then
\[\mu(L_1,L_2) = \mu(K_1,K_2).\]

\end{enumerate}

Proofs of these facts follow from the interpretation of the Maslov index as an
intersection number and can be found in \cite{CLM4}. It is worth noting that  the
proof of our main theorem  requires only these elementary properties of the Maslov
index and avoid more technical tools such as symplectic reduction.   

Path additivity does not hold with all possible conventions, and it is the reason
we use the chosen convention. There are other conventions; a popular choice is to
use the $(-\epsilon,\epsilon)$ Maslov index  and spectral flow, since this choice
fits in with the index theorem of Atiyah-Patodi-Singer as stated in
\cite{APS}.  To go back and forth between conventions one needs only to know that
if $\mu'$ is another convention, then there exist numbers $\sigma_0$ and
$\sigma_1$ in $\{-1,0,1\}$ and $e\in \{1,-1\}$ so that
$$\mu'(L,M)=e\cdot\mu(L,M)+\sigma_0\cdot\dim(L(0)\cap M(0)) +
\sigma_1\cdot\dim(L(1)\cap M(1)).$$  A similar remark applies to the spectral
flow, and it is not hard to see that  the formula of our main result, Theorem
\ref{splitthm}, remains true provided one chooses the spectral flow and Maslov
index conventions compatibly,   after perhaps adding a correction term depending
only on the dimensions of $\ker D(0)$ and $\ker D(1)$.

 Two further simple facts we will use in Section \ref{UsersGuide} without
explicit mention are that (with our chosen conventions), 

\begin{enumerate}
\item $\mu(L,M)=0$ if $L$ and $M$ are constant paths.
\item  If $L,M$ are paths of Lagrangians in $H_{\nu}(S)$, then 
$$\mu(P^-_\nu\oplus L, M\oplus P^+_\nu)=\mu(L,M).$$
\end{enumerate} These are easy consequences of the definitions.

The following remarkable theorem of Nicolaescu will be the basis of what follows.

\begin{theorem}\label{NicolaescusSplittingTheorem} Let $D(t),\ t\in [0,1]$ be a
smooth path of (cylindrical, neck-compatible, self-adjoint) Dirac operators on a
smooth, oriented, closed, odd dimensional Riemannian manifold $M$ which splits as
$M=X\cup_{\Sigma}Y$. Then
\[\SF(D) = \mu(\Lambda_X(D),\Lambda_Y(D)).\]
\end{theorem}

  So the theorem explicitly states the intuitively appealing idea that counting
kernel elements along the path ({\em i.e.} counting eigenvalues that cross through
zero) is equivalent to counting pairs of boundary values that match up ({\em i.e.}
nontrivial intersections between the Cauchy data spaces). Theorem
\ref{NicolaescusSplittingTheorem} was first proved by Nicolaescu for paths of
Dirac operators whose endpoints have trivial kernel in \cite{Nico}.  The
restriction to trivial kernel at the endpoints was removed in
\cite{AMD2}.

A similar theorem may be stated for manifolds with boundary.  In this case we must
impose boundary conditions for the spectral flow to be well defined.   This is the
subject of the next definition.  
\begin{definition}  Let $X$ be a manifold with boundary $\partial X=\Si$ and $D$ a
self-adjoint Dirac operator on $X$ in cylindrical form with tangential operator
$S$. A {\it self-adjoint elliptic boundary condition} is a Lagrangian subspace
$B\subset L^2(E_{|\Si})$ which contains $P^+_\nu(S)$ as a finite codimensional
subspace for some
$\nu$.  \end{definition} See \cite{BW} and \cite{Nico2} for details.  The
condition that
$B$ be  Lagrangian implies that the operator $D$ on $X$ with boundary conditions
$B$ is self-adjoint.  The requirement that $B$ contain $P^+$ with finite
codimension ensures that the operator $D$ acting of sections over $X$ whose
restriction to the boundary lies in $B$ is elliptic.  Thus given a path $D(t)$ of
Dirac operators on $X$ and a path of elliptic self-adjoint boundary conditions
$B(t)$ the spectral flow $SF(D,B)$ is defined.

Then Nicolaescu's theorem extends to the bounded case as follows.

\begin{theorem}\label{NicolaescusSplittingTheoremWithBoundary} Let $D(t),\ t\in
[0,1]$ be a smooth path of (cylindrical, neck-compatible, self-adjoint) Dirac
operators  on a smooth, oriented, odd dimensional Riemannian manifold $X$ with
nontrivial boundary $\partial X=\Si$. Let
$B(t)$ be a smooth path of elliptic boundary conditions for $D(t)$.  Then
\[\SF(D,B) = \mu(\Lambda_X(D),B).\]
\end{theorem}

\section{The general splitting formula} In this section we state and prove the
general splitting formula.  The formula expresses the spectral flow of a path of
Dirac operators on a closed manifold in terms of the spectral flows of the
restricted paths (with associated elliptic boundary conditions).  Whereas other
results of this type have many additional hypotheses and produce more succinct
formulas, our result requires only the minimal hypotheses, but produces a longer
formula.  In Section
\ref{UsersGuide} we discuss additional conditions that may be imposed to make
various terms in our formula vanish or cancel.

The set-up is as follows.  Let $D(t)$ be a smooth path of Dirac operators on a
smooth, oriented, closed, odd dimensional Riemannian manifold $M$.  Suppose that
$M$ can be split along a hypersurface $\Sigma$ ($M=X\cup_{\Sigma}Y$) and that each
$D(t)$ is cylindrical and neck-compatible with respect to this splitting. Let
$B_X(t)$ and $B_Y(t)$ be paths of elliptic boundary conditions for
$D(t)$ restricted to $X$ and $Y$ respectively. Then we will show that  there is an
11 term formula

\begin{equation}\label{firstformula} SF(D) = SF(D_X,B_X) + SF(D_Y,B_Y) +
\mu(B_Y(1-t),B_X(1-t)) + \sum_{i=1,2,4,5,7,8,10,11}\mu(L_i,M_i)
\end{equation} The   $\mu(L_i,M_i)$  are certain Maslov indices. They will be
defined  below and will be discussed at length in the following section.  

\vskip4ex

 Theorem \ref{NicolaescusSplittingTheorem}  allows us to replace 
  $SF(D)$ by
$\mu(\Lambda_X(D),\Lambda_Y(D))$.  We have at our disposal the path additivity and
the homotopy invariance of the Maslov index. We will describe paths $L$ and $M$
that are homotopic rel endpoints to
$\Lambda_X(D)$ and
$\Lambda_Y(D)$ respectively.  These new paths will each be the concatenation of
eleven pieces ($L_i$ and $M_i$ respectively).  Each piece will contribute a term 
to the right hand side of Formula \ref{firstformula}.

To begin let $\nu_0\ge 0$ and $\nu_1\ge 0$ be numbers chosen so that 
\begin{enumerate}
 \item $\nu_0$ is  in the nonresonance range for $D(0)$ on $X$
and the tangential operator $S(0)$ has a spectral gap at $\nu_0$, and
\item $\nu_1$ is in the nonresonance range for $D(1)$ on $Y$ and
the tangential operator $S(1)$ has a spectral gap at $\nu_1$.
\end{enumerate}

 We abbreviate the notation for the Cauchy data spaces using the symbol
$\Lambda_X^r(t)$ for
$\Lambda_{X^r}(D(t))$.   Moreover $\La_X(t)$ means
$\La_{X^0}(D(t))=\La_{X}(D(t))$. Similar notation applies to $Y$.    Nicolaescu's
adiabatic limit theorem (Theorem
\ref{NicoAdiabaticLimitTheorem} above) shows that there exists a  Lagrangian 
$L_X(0)\subset H_{\nu_0}$ (and gives a recipe for constructing it) so that
$$\lim_{r\to \infty} \La^r_X(0)=P^-_{\nu_0}(S(0))\oplus L_X(0)$$ and  there
exists a Lagrangian $L_Y(1)\subset H_{\nu_0}$ so that
$$\lim_{r\to \infty} \La^r_Y(1)= L_Y(1)\oplus P^+_{\nu_1}(S(1)).$$

We can now enumerate the eleven pieces of each path. 
\begin{enumerate}
\item Let $L_1$ be the path  starting at
$\Lambda_{X}^0(0)$  and ending at 
$\lim_{r\to\infty}\Lambda_{X}^r(0)=P^-_{\nu_0}(S(0))\oplus L_X(0) $ obtained by
stretching.  An explicit formula is given in the statement of Lemma
\ref{stretchpath}.       Let $M_1$ be the constant path at $\Lambda_Y(0)$.

\item  Let $L_2$ be any path of Lagrangians starting at $P^-_{\nu_0}(S(0))\oplus
L_X(0)$ and ending at $B_Y(0)$ so that for all $t$,   $L_2(t)$ is a self-adjoint
elliptic boundary condition for the restriction of $D(0)$ to $Y$ (or more
generally  it suffices to assume that
$(L_2(t), \La_Y(0))$ are a Fredholm pair).   Let  $M_2$ be the constant path
$\Lambda_Y(0)$.

\item Let $L_3(t) $ be  $B_Y(t)$ and let   $M_3(t)$ be $\Lambda_Y(t)$.  Theorem
\ref{NicolaescusSplittingTheoremWithBoundary} applied to $Y$ implies that 
\begin{equation}\mu( L_3,M_3)=SF(D_{|Y},B_Y).
\end{equation}

\item Take $L_4$ to be the constant path  $B_Y(1)$, and  let $M_4$ be the path from
$\Lambda_Y(1)$ to $\lim_{r\to\infty}\Lambda_Y^{r}(1)= L_Y(1)\oplus
P^+_{\nu_1}(S(1))$ obtained by stretching as in Lemma \ref{stretchpath}.

\item  Let $L_5$ be the constant path $B_Y(1)$.  For $M_5$ choose a path of
Lagrangians  starting at $L_Y(1)\oplus P^+_{\nu_1}(S(1))$  and ending at  $B_X(1)$
so that for all $t$,
$M_5(t)$ is a self-adjoint elliptic boundary condition for the restriction of
$D(1)$ to
$X$ (or more generally so that $(\La_X(1), M_5(t))$ form a Fredholm pair).

\item  Let $L_6$ be the   path as $L_3$ run backwards, i.e. $L_6(t)=L_3(1-t)$, and
let 
$M_6$ be
$B_X$ run backwards.   Thus 
\begin{equation}\mu(L_6,M_6)=\mu(B_Y(1-t),B_X(1-t))\end{equation}

\item Let $L_7$ be $L_2$ run backwards and let $M_7$ be the constant path
$B_X(0)$.  
\item Let $L_8$ be $L_1$ run backwards and $M_8$   the constant path
$B_X(0)$. 
\item Let $L_9$ be the path $\Lambda_X(t)$ and   let $M_9$ be the path $B_X(t)$.  
Theorem 
\ref{NicolaescusSplittingTheoremWithBoundary} says: 
\begin{equation}   \mu(L_9,M_9)=SF(D_{|X},B_X).
\end{equation} 

\item Take $L_{10}$ to be the constant path $\Lambda_X(1)$ and  $M_{10}$ to be
$M_5$ run backwards.  
\item Finally, let $L_{11}$ be the constant path $\Lambda_X(1)$ and  $M_{11}$ to
be  $M_4$ run backwards.
\end{enumerate}

The reader may verify that the composite path $L=L_1*L_2*\cdots *L_{11}$ is
defined and is homotopic rel endpoints to the path $\La_X$.   Similarly
$M=M_1*M_2*\cdots * M_{11}$ is  homotopic rel endpoints to $\Lambda_Y$.  Hence 
$$SF(D)=\mu(\La_X,\La_Y)=\mu(L,M)=\sum_{i=1}^{11}\mu(L_i,M_i),$$ using  homotopy
invariance of the Maslov index and  additivity of the Maslov index under
composition of paths.

We summarize our conclusions in the following theorem.

\begin{theorem}\label{splitthm}  Let $D(t)$ be a continuous path of self-adjoint
Dirac operators on a smooth, closed, oriented, odd dimensional Riemannian manifold
$M$.  Suppose that $M$ can be split along a hypersurface $\Si$ ($M=X\cup_{\Si}Y$)
and that each
$D(t)$ is cylindrical and neck compatible with respect to this splitting.  Let
$B_X(t)$ and
$B_Y(t)$ be paths of self-adjoint elliptic boundary conditions for the restriction
of
$D(t)$ to $X$ and $Y$.

Then 
\[ SF(D)= SF(D_{|X},B_X) + SF(D_{|Y},B_Y) +\mu(B_Y(1-t),B_X(1-t)) +\sum_{i\ne
3,6,9}\mu(L_i,M_i).\]
\end{theorem} 

\section{User's guide to Theorem \ref{splitthm} }\label{UsersGuide}

 In this section we explain how to use Theorem \ref{splitthm}.  Specifically we
show how various natural hypotheses simplify the formula, and then derive some
earlier theorems as  consequences. We will not exhaust all the possibilities, but
hope to give some indication of the utility of the formula.  

The authors' background  concerns the application of this subject to the odd
signature operator coupled to a path of connections starting and ending at flat
connections. This is the kind of operator considered in topological applications
of spectral flow, such as computations of Atiyah-Patodi-Singer $\rho_\alpha$
invariants, Casson's invariant, and Floer homology.  The methods we describe are
particularly well suited for this class of problem.

\subsection{Transversality at endpoints and stretching}

First some notation.  We have defined $X^r$ and $Y^r$ to be the manifolds
obtained by adding a collar of length $r$ to $X$ and $Y$.  Let $M^r$ be the
closed manifold obtained by stretching $M$ along $\Si$, so 
$$M^r=X^r\cup_{\Si} Y^r.$$

\begin{hypothesis}\label{trans} The adiabatic limits of the Cauchy data spaces are
transverse at the endpoints.
$$\lim_{r\to \infty}\La_X^r(i)\cap \lim_{r\to \infty}\La^r_Y(i)=0, \ \ i=0,1.$$
\end{hypothesis}

\begin{proposition} \label{sttt} Suppose that Hypothesis \ref{trans} holds.  Then
there exists an $r_0\ge 0$ so that replacing $M$ by $M^r$  for $r\ge r_0$ in
Theorem
\ref{splitthm}, the terms $\mu(L_1,M_1)$ and $\mu(L_{11}, M_{11})$ vanish.
\end{proposition}
\proof   Continuity of the path of Lemma \ref{stretchpath} implies that  there
exists some $r_0$ so that  the Lagrangians $\La^r_X(i)$ and $\La^r_Y(i)$ are
transverse for $r\ge r_0$ and
$i=0,1$.   Then the Lagrangians $L_1(t)$ and $M_1(t) $ are transverse for all
$t\in [0,1]$ and hence $\mu(L_1,M_1)=0$.  The same argument applies at the other
end of the path to show that $\mu(L_{11},M_{11})=0$. \qed

Notice that the two cases are independent, i.e. if the limits of the Cauchy data
spaces are transverse at the initial point then $\mu(L_1,M_1)=0$ for $r$ large
enough, and if they are transverse at the terminal point then
$\mu(L_{11},M_{11})=0$ for
$r$ large enough.

A slight generalization of this can be obtained by using the following hypothesis.

\begin{hypothesis}\label{const} For $i=0$ and $1$, the dimension of
$\La_X^r(i)\cap
\La_Y^r(i)$ is independent of $r$  for $r\ge r_0$ and equals  the dimension of
the intersection  of the limits of the Cauchy data spaces
$$\dim(\La_X^r(i)\cap \La_Y^r(i))=\dim(\lim_{r\to\infty}\La_X^r(i)\cap
\lim_{r\to\infty}\La_Y^r(i)).$$
\end{hypothesis}

Notice that the intersection $\La^r_X(i)\cap \La^r_Y(i)$ is isomorphic to the
kernel of
$D(i)$ on $M^r$, so Hypothesis \ref{const} implies  (but is in general stronger)
that  the dimension of this kernel is independent of $r$.

\begin{proposition} If Hypothesis \ref{const}  holds then  after  replacing $M$ by
$M^r$  for
$r\ge r_0$ in Theorem \ref{splitthm}, the terms $\mu(L_1,M_1)$ and $\mu(L_{11},
M_{11})$ vanish.
\end{proposition}
\begin{proof} Let $\La^{\infty}_X(0)$ denote the adiabatic limit of $\La_X^r(0)$
with similar notation for $Y$.   

Fix $r\ge r_0$ and let  $u\ge r$. 
Since 
 $\dim(\La^u_X(0)\cap\La^{r}_Y(0))$ is isomorphic to the kernel of $D$ on
$M^{u+r}$, which in turn is isomorphic to
$\dim(\La^{(u+r)/2}_X(0)\cap\La^{(u+r)/2}_Y(0))$, Hypothesis
\ref{const} implies that  
\[\dim(\La^{u}_X(0)\cap\La^{r}_Y(0)) = \dim(\La^{\infty}_X(0)\cap
\La^{\infty}_Y(0)).\] 
Thus the dimension of the intersection of $L_1(t)$ with $M_1(t)$ is independent
of
$t$.  This implies that $\mu(L_1,M_1)=0$.  A similar argument shows that
$\mu(L_{11},M_{11})$ vanishes.\end{proof}

\subsection{Choice of Boundary Conditions}

The boundary conditions $B_X$ and $B_Y$ can be restricted to simplify the splitting
formula.  The most direct way to do this is just to kill the terms $\mu(L_2,M_2)$,
$\mu(L_5,M_5)$, $\mu(L_7,M_7)$, and $\mu(L_{10},M_{10})$ by choosing the boundary
conditions $B_Y(0)$ and $B_X(1)$ as follows.

\begin{hypothesis}\label{kill}
$B_Y(0)=P^-_{\nu_0}(0)\oplus L_X(0)$ and $B_X(1)=L_Y(1)\oplus P^+_{\nu_1}(1)$.
\end{hypothesis}

\begin{proposition} \label{kill2} Assume that Hypothesis \ref{kill} holds. Then
one can choose the paths $L_2$ and $M_5$ (and  their reverses $L_7$ and
$M_{10}$)  so that
$$\mu(L_2,M_2)=\mu(L_5,M_5)=\mu(L_7,M_7)=\mu(L_{10},M_{10})=0.$$
\end{proposition}
\proof  Take $L_2$ and $M_5$ to be constant paths.  Then $L_7$ and $M_{10}$ are
also constant.  By definition, $M_2$, $L_5$, $M_7$, and $L_{10}$ are constant. 
Thus the four terms are Maslov indices of constant paths, and so all vanish.\qed 

We could  have taken the point of view in Theorem \ref{splitthm} that only
boundary conditions satisfying  Hypothesis \ref{kill} are allowed. This would
have given a formula with four fewer terms, but the result would have been less
flexible. The decision to state the theorem as we did was made to decouple the
choice of boundary conditions from the analysis of the adiabatic limits of the
Cauchy data spaces.

\subsection{The nonresonance range, limiting values of extended $L^2$ solutions,
and adiabatic limits}

We next give a slightly more detailed description of the adiabatic limit
$\lim_{r\to
\infty}\La_X^r$ which can be useful in controlling some of the terms.

\begin{definition}\label{extended} Let $D$ be a cylindrical Dirac operator as
above on a manifold $X$ with boundary.  The Lagrangian subspace
$$\wtilde{L}_X(D)\subset\ker S$$ defined to be the symplectic reduction of the
Cauchy data space to the kernel of $S$
 $$\wtilde{L}_X(D)=\proj_{\ker S}\left(\La_X(D) \cap(\ker S \oplus
P^+(S))\right)$$ is called the {\it limiting values of extended $L^2$ solutions.}
(This terminology comes from
\cite{APS}.)\end{definition}

For convenience, we recall the notation for several Lagrangians that appear in
this section.
\begin{enumerate}
\item $\La_X^r$, the Cauchy data space on $X^r$. This is an infinite dimensional
Lagrangian subspace of $L^2(E|_{\Si})$.
 \item $\wtilde{\La}_X$, the symplectic reduction of the (length $0$) Cauchy data
space
$\La_X^0$ to
$H_{\nu_0}$ (Equation
\ref{eq1}) where $\nu_0\ge 0$ is greater than or equal to the nonresonance level
of $D$ and $S$ has a spectral gap at $\nu_0$.   This is a finite dimensional
Lagrangian subspace of the  symplectic vector space
$H_{\nu_0}$ defined in Equation \ref{hnu}.
\item $L_X$, the limit of $e^{-rS}\wtilde{\La}_X$ as $r\to\infty$, a Lagrangian
subspace of
$H_{\nu_0}$ (Equation \ref{defofL}). Thus the adiabatic limit
$\lim_{r\to\infty}\La_X^r=P^-_{\nu_0}\oplus L_X$.
\item $\wtilde{L}_X$, the limiting values of extended $L^2$ solutions, defined as
the symplectic reduction of the Cauchy data space $\La_X^0$ to the kernel  of $S$
in Definition
\ref{extended}.
\end{enumerate}

The following Theorem relates these Lagrangians, and indicates the
structure of
$L_X$.  It is convenient to extend the notation slightly, so that for the
statement and proof of this Theorem we will  allow $\nu<0$ in the definition of
$P^+_\nu$ (Equation
\ref{pplusnu}). For example if  $\nu$ is positive and in the complement of the
spectrum of
$S$,   then $H_{\nu}\oplus P_{\nu}^+= P^+_{-\nu}.$

Notice that there is a descending filtration of $H_{\nu_0}\oplus P_{\nu_0}^+$
corresponding to the increasing list of eigenvalues
$-\la_{n+1}<-\la_n\cdots<0<\la_1<\cdots$ assuming that $\la_n\le \nu_0<\la_{n+1}$.
\begin{equation}\label{filt} P^+_{-\lambda_{n+1}}\supset
P^+_{-\lambda_{n}}\supset\cdots\supset P^+_0\supset P^+_{\la_1}\supset\cdots
\end{equation}

\begin{theorem}\label{tricky} With notation as described above, the following
statements are true.
\begin{enumerate}
\item  $L_X=\lim_{r\to\infty}e^{-rS}\wtilde{\La}_X$.
\item   $\wtilde{L}_X\subset L_X$. 
\item If $\nu_0=0$, then $L_X=\wtilde{L}_X=\wtilde{\La}_X$.
\item Let $0<\la_1<\la_2<\cdots<\la_n<\la_{n+1}$ denote the ordered list of 
positive eigenvalues of the tangential operator $S$  so that
$\la_n\le\nu_0<\la_{n+1}$. Let
$E_i^\pm $ denote the   $\pm \la_i$ eigenspace. 

 Then 
\begin{equation}\label{decompo}  H_{\nu_0}=E_n^-\oplus E_{n-1}^-\oplus\cdots\oplus
\ker S
\oplus E^+_1\oplus\cdots
\oplus E^+_n\end{equation}
 and the Lagrangian $L_X$ decomposes in this direct sum in the form
\begin{equation}\label{ass}L_X= W_n\oplus  W_{n-1}\oplus \cdots\wtilde{L}_X\oplus
JV_1\oplus \cdots\oplus JV_n\end{equation}
 where
$W_i\subset E_i^-$ are subspaces, $V_i\subset E_i^-$ are their orthogonal
complement in
$E_i^-$ (and so $JV_i\subset E_i^+$).  Moreover this decomposition  exhibits  
$L_X$ as  the associated graded to the  filtration  of $\La_X\cap (H_{\nu_0}\oplus
P^+_{\nu_0})=\La_X\cap P^+_{-\la_{n+1}}$ obtained by intersecting  $\La_X$  with
the decreasing filtration
  given in Equation \ref{filt}.
\end{enumerate}

\end{theorem}
\proof The first assertion is the definition of $L_X$.   For the third assertion,
if $\nu_0=0$, then $\wtilde{L}_X=\wtilde{\La}_X$ by definition.  Since the
operator $S$ is zero on its kernel, the restriction of 
$e^{-rS}$ to $\ker S$ is the identity so that $L_X=\wtilde{\La}_X$.  The second
assertion follows from the fourth.  So we will prove the fourth.

This result follows from a more careful analysis of the flow to the adiabatic
limit.  Notice that because the $P^+_\nu$ are defined in terms of strict
inequalities, 
$P^+_{-\la_i}$ is the span of the eigenvectors whose eigenvalues are {\it greater}
than
$-\la_i$. Thus $E_i^-=P^+_{-\la_{i+1}}/P^+_{-\la_{i}}$ and
$E_i^+=P^+_{\la_{i-1}}/P^+_{\la_{i}}$

Let 
$$W_n=\proj_{E_n^-}(\La_X\cap P^+_{-\lambda_{n+1}})\subset E_n^-.$$  Thus
$$W_n= (\La_X\cap P^+_{-\lambda_{n+1}})/(\La_X\cap P^+_{-\lambda_{n}}).$$ Next let
$$W_{n-1}=\proj_{E_{n-1}^-}(\La_X\cap P^+_{-\lambda_{n}})\subset E_{n-1}^-.$$
Continue, peeling one space off at a time in the decomposition  \ref{decompo} of
$H_{\nu_0}$. We change notation when we get to $\ker S$ to be consistent with our
previous notation,  Thus (by definition)
$$\wtilde{L}_X=\proj_{\ker S}(\La_X\cap P^+_{-\lambda_{1}})\subset \ker S.$$
Continue by letting 
$$V_1'=\proj_{E_1^+}(\La_X\cap P^+_{0})\subset E_1^+,$$
$$V_2'=\proj_{E_2^+}(\La_X\cap P^+_{\lambda_{1}})\subset E_2^+,$$ and so forth
until the last step
$$V_n'=\proj_{E_n^+}(\La_X\cap P^+_{\lambda_{n-1}})\subset E_n^+.$$

Suppose that $(w_n, w_{n-1},\cdots, w_1,h,v_1,\cdots,v_n, q)$ is an element of
$\La_X\cap(H_{\nu_0}\oplus P^+_{\nu_0})$ expressed in the decomposition 
\ref{decompo} with the additional element  $q\in P^+_{\nu_0}$.  Then
$(w_n, w_{n-1},\cdots, w_1,h,v_1,\cdots,v_n)$ is in  $
\wtilde{\La}_X$.  

Either $w_n=0$ or else $w_n\in W_n-\{0\}$.  Since
$${1\over e^{r\lambda_n}}(w_n, w_{n-1},\cdots, w_1,h,v_1,\cdots,v_n)\in
\wtilde{\La}_X$$ and since $e^{-rS}$ acts on the decomposition  of
Equation \ref{decompo} diagonally with (decreasing)  eigenvalues
$e^{r\lambda_n},e^{r\lambda_{n-1}},\cdots$ it follows that  if $  w_n\ne 0$
$$\lim_{r\to\infty}e^{-rS}{1\over e^{r\lambda_n}}(w_n, w_{n-1},\cdots,
w_1,h,v_1,\cdots,v_n) =(w_n,0,0,\cdots,0)$$ and so $w_n$ is in $\lim_{r\to
\infty}e^{-rS}\wtilde{\La}_X= L_X.$ Arguing by induction one obtains 
$$L_X=W_n\oplus\cdots\oplus W_1\oplus \wtilde{L}_X\oplus V_1'\oplus\cdots\oplus
V_n'.$$

We must now see that $V_i'=JV_i$ where $V_i$ is the orthogonal complement of $W_i$
in
$E_i^-$.  But this follows from the fact that $L_X$ is a Lagrangian subspace and
dimension counting.  Indeed, since the symplectic structure on
$H_{\nu_0}$ is given by $\omega(x,y)=<x,Jy>$,  $JV_i'$ is orthogonal to $W_i$ and
so lies in
$V_i$. If for some $i$, $JV_i'$ were a proper subspace of $V_i$ then by counting
dimensions (and using the fact that the limiting values of extended
$L^2$ solutions $\wtilde{L}_X\subset
\ker S$ is a Lagrangian subspace of $\ker S$) it would follow that $L_X$ has too
small a dimension to be a Lagrangian. Thus $JV_i'=V_i$ and so $V_i'=JV_i$ as
claimed.

The assertion that  $L_X$ is the associated graded to the filtration is simply a
brief description of  how the
$W_i,\wtilde{L}_X, V_i'$ were constructed.  
\qed

The fourth statement of  Theorem \ref{tricky} suggests a more useful and
sophisticated alternative to Hypothesis
\ref{kill}.   The underlying motivation comes from the fact that it is much easier
to calculate $\wtilde{L}_X\subset \ker S$ than to calculate $L_X\subset
H_{\nu_0}$. Even getting a handle on the  nonresonance level $\nu_0$  can be a
difficult problem. 

There is a natural choice of path of Lagrangians starting at $P^-_{\nu_0}\oplus
L_X$ and ending at $P^-\oplus\wtilde{L}_X$ defined  as follows.   Notice that the
symplectic subspaces $E_i^-\oplus E_i^+$ have a further decomposition as a direct
sum
$$E^-_i\oplus E_i^+ = (W_i\oplus V_i)\oplus (JW_i\oplus JV_i) = (W_i\oplus
JW_i)\oplus (V_i\oplus JV_i).$$

Use these decompositions to define a path $C(t)$ by the formula
\begin{equation}\label{rot} C:t\mapsto P^-_{\nu_0}\oplus\wtilde{L}_X\oplus (W_1
\oplus e^{-(1-t){\pi\over 2}  J}V_1)\oplus\cdots
\oplus(W_n \oplus e^{-(1-t){\pi\over 2}  J}V_n)
\end{equation} 

Then we can make the hypothesis
\begin{hypothesis} \label{confusing}
$B_Y(0)=P^-\oplus\wtilde{L}_X$, $B_X(0)= A\oplus P^+$ for some Lagrangian
$A\subset \ker S(0)$,  and
$L_2(t)=C(t)$.  
\end{hypothesis}

\begin{lemma}\label{cann} If Hypothesis \ref{confusing} holds, then
$\mu(L_7,M_7)=0$.
\end{lemma} 
\proof  This is a consequence of the conventions  we are using for spectral flow. 
Recall that
$L_7$ is $L_2$ run backwards, and $M_7$ is the constant path at $B_X(0)$.

Using Theorem \ref{tricky} we see that the Lagrangians
$L_7(t)=C(1-t)$ and $M_7(t)=B_X(0)$ intersect in the direct sum
$$L_7(t)\cap M_7(t)=\begin{cases} (\wtilde{L}_X\cap A) \oplus JV_1\oplus
JV_2\oplus\cdots\oplus JV_n&\hbox{ if } t=1\\
\wtilde{L}_X\cap A&\hbox{ if } 0\le t < 1.\end{cases}$$

Thus $e^{\epsilon J}L_7(t)$ is transverse to $M_7(t)$ for all $t$, so that
$\mu(L_7,M_7)=0$.\qed

A similar argument shows that if $B_X(1)=\wtilde{L}_Y\oplus P^+$, $B_Y(1)=
P^-\oplus A$  for some Lagrangian $A\subset \ker S(1)$,  and
$M_5$ is chosen in a manner similar to $C(t)$ above, then $\mu(L_5,M_5)=0$.   

Finally, with these choices and some additional transversality conditions, one can
sometimes also  compute $\mu(L_2,M_2)$ and $\mu(L_{10}, M_{10})$ in terms of the
sum of the dimensions   of the $V_i$ (which is the same as the dimension of the
$L^2$-kernel of $D$ on $X^\infty$; see below) after a preliminary stretching. 
Lemma \ref{cann} will be used in a slightly different context in Theorem
\ref{wbound} below.

\subsection{When the $L^2$ kernel of $D_{|X}$  or $D_{|Y}$ vanishes at the
endpoints}

The nonresonance level for $D_{|X}$ is zero if $\La_X\cap P^+=0$.  This is
equivalent (see
\cite{APS}) to the vanishing of  the $L^2$-kernel of the natural extension  of
$D_{|X}$ to $X^\infty=X\cup \Si\times [0,\infty)$.

Nicolaescu's adiabatic limit theorem in this context says   that   if $\La_X\cap
P^+=0$ (i.e.  if the $L^2$ kernel of $D_{|X}$ on $X^\infty$ is zero), then
$\lim_{r\to
\infty}\La^r_X = P^-\oplus \wtilde{L}_X$

\begin{hypothesis}\label{L2ker} The operators $D(i)_{|X^\infty}$ and
$D(i)_{|Y^\infty}$ have no $L^2$ kernels for $i=0$ and $1$.
\end{hypothesis} (In the  terminology of \cite{Nico} the operators
$D(i)_{|X^\infty}$ and
$D(i)_{|Y^\infty}$ are non-resonant.)

If Hypothesis \ref{L2ker} holds, then Hypothesis \ref{trans} holds if and only if
$ \wtilde{L}_X(0)\cap  \wtilde{L}_Y(0)=0$ and $  \wtilde{L}_X(1)\cap 
\wtilde{L}_Y(1)=0$.  Moreover, in this case  one can satisfy Hypothesis \ref{kill}
by  letting
$B_Y(0)=P^-(0)\oplus\wtilde{L}_X(0)$ and $B_X(1)=\wtilde{L}_Y(1)\oplus P^+(1)$.

Let us use these ideas to give a simple proof  of a theorem of Bunke \cite{bunke}
(see also
\cite{CLM2}, Theorem A).

Consider the case  when the tangential operator has no kernel along the path.  The
following theorem appears (in  different notation) in \cite{bunke}, it also
follows from Theorem A of \cite{CLM2}.

\begin{theorem}\label{bunke} Suppose that the kernel of the tangential operator
$S(t)$ vanishes for all
$t$. Suppose that Hypothesis \ref{L2ker} holds (at the endpoints).

Then there exists an $r_0$ so that for $r\ge r_0$, 
$$SF(D, M^r)=SF(D_{X^r}; P^+) + SF(D_{Y^r}; P^-).$$
\end{theorem}

\proof Hypothesis \ref{L2ker}  together with the vanishing of the kernels of the
tangential operators
$S(0)$ and $S(1)$ imply that the adiabatic limits are
$$\lim_{r\to \infty}\La^r_X(i)=P^-(i)$$ and 
$$\lim_{r\to \infty}\La^r_Y(i)=P^+(i)$$ for $i=0$ and $1$. Since $P^+(i)$ is
transverse to $P^-(i)$,  Hypothesis \ref{trans} holds.  Thus, Proposition
\ref{sttt}  implies that $\mu(L_1,M_1)$  and
$\mu(L_{11},M_{11})$ vanish after sufficient stretching.

Since the kernel of $S(t)$ is zero for all $t$, the spaces $P^\pm(t)$ vary
continuously (\cite{KK-illinois, BW}) and so we can take $B_X(t)= P^+(t)$ and
$B_Y(t)=P^-(t)$. This immediately implies that $\mu(L_i,M_i)$ vanishes for $i=2,5,
7$ and $10$ according to Proposition \ref{kill2} since Hypothesis \ref{kill}
holds.  

We also have $\mu(L_6,M_6)=0$ since this equals $\mu(P^-(1-t),P^+(1-t))$ and 
$P^+(t)$ is transverse to $P^-(t)$ for all $t$. 

The path $L_4$ is the constant path at  $B_Y(1)=P^-(1)$.  The path $M_4$ is is the
path from $\La_Y(1)$ to $\lim_{r\to\infty} \La_Y^r(1)=P^+(1)$.  Since $P^-(1)$ is
transverse to $P^+(1)$,   after perhaps making
$r$ larger, $\mu(L_4,M_4)=0$. Similarly, after perhaps making $r$ larger,
$\mu(L_8,M_8)=0$. 

The only terms remaining are $\mu(L_3,M_3)= SF(D_{|Y^r};P^-)$ and
$\mu(L_9,M_9)=\SF(D_{|X^r}, P^+)$. This completes the proof.\qed

We next state and prove the theorem of Yoshida and Nicolaescu. The proof we give is
identical to  Nicolaescu's. We include it for the convenience of the reader and
to introduce another useful technique which can be combined with the our
methods.   Theorem
\ref{3term} below generalizes both Theorem \ref{bunke} and Theorem \ref{yoshida}.

Notice that if the operator 
$D(t)_{|X^\infty}$ has no
$L^2$ kernel for all
$t\in[0,1]$,  and if the kernel of the tangential operator $S(t)$ has constant
dimension  along the path, then the limiting values of extended $L^2$ solutions 
 $\wtilde{L}_X(t)$  is a continuous path of  Lagrangians.  (The symplectic
reduction is {\it clean} in the sense of \cite{Nico}.)  This is because for all
$t$ the projection
$$\La_X(t)\cap(\ker S(t)\oplus P^+(S(t)) )\to \ker S(t)$$ (with image
$\wtilde{L}_X(t)$) has no kernel and it is easily checked that the image is
continuous since the path
$\La_X(t)$ is continuous.

Thus the statement of the following theorem makes sense.  This is Corollary 4.4 in
\cite{Nico}.

\begin{theorem} \label{yoshida} (Yoshida, Nicolaescu)  Suppose that for all
parameters
$t\in [0,1]$, the operators $D(t)_{|X^\infty}$ and $D(t)_{|Y^\infty}$ have no
$L^2$ kernel.  Assume furthermore that the  kernel of the tangential operator
$\ker S(t)$ has constant dimension, i.e. independent of $t\in[0,1]$. Assume that 
\begin{equation}\label{eqqq1}\wtilde{L}_X(0)\cap \wtilde{L}_Y(0)=0\end{equation}
 and 
\begin{equation}\label{eqqq2}\wtilde{L}_X(1)\cap \wtilde{L}_Y(1)=0.\end{equation}

Then there exists an $r_0\ge 0$ so that for $r\ge r_0$,
$$SF(D, M^r)= \mu(\wtilde{L}_X,\wtilde{L}_Y).$$
\end{theorem}

\proof  Since $0$ is the nonresonance level,
$\lim_{r\to\infty}\La_X^r(t)=P^-(t)\oplus \wtilde{L}_X(t)$ and
$\lim_{r\to\infty}\La_Y^r(t)=\wtilde{L}_Y(t)\oplus P^+(t)$ for all $t\in[0,1]$.
Together with   Equations
\ref{eqqq1} and
\ref{eqqq2} this implies  that if $r$ is large enough, $\La_X^r(i)$ is transverse
to
$\La_Y^r(i)$ for $i=0$ and $1$. Fix $r_0\ge 0$ so that they are transverse  for
all
$r\ge r_0$. let $r_1$ be any number greater than or equal to $r_0$.

Consider the homotopy from the path of pairs $(\La_X^{r_1}(t),\La_Y^{r_1}(t))$
to 
$(P^-(t)\oplus \wtilde{L}_X(t),\wtilde{L}_Y(t)\oplus P^+(t))$ obtained by letting
$r$ go to infinity. This is not a rel endpoints homotopy, but does exhibit a rel
endpoint homotopy from the path of pairs $(\La_X^{r_1}(t),\La_Y^{r_1}(t))$ to the
composite of three paths:
$$A(t)= \begin{cases} (\La^{r_1/1-t}_X(0), \La^{r_1/1-t}_Y(0))     & \hbox{ if }
t<1,\\
                (P^-(0)\oplus \wtilde{L}_X(0), \wtilde{L}_Y(0)\oplus P^+(0))
 & \hbox{ if } t=1,\end{cases}$$
$$B(t)=(P^-(t)\oplus \wtilde{L}_X(t),\wtilde{L}_Y(t)\oplus P^+(t)),$$ and
$$C(t)=\begin{cases}  
        (P^-(1)\oplus \wtilde{L}_X(1), \wtilde{L}_Y(1)\oplus P^+(1)) & \hbox{ if }
t=0,\\
         (\La^{r_1/t}_X(1),\La^{r_1/t}_Y(1))&\hbox{ if } t>0,
\end{cases}$$ and so $SF(D, M^{r_1})=\mu(A*B*C)=\mu(A)+\mu(B)+\mu(C)$.

Since $\La_X^r(i)$ is transverse to
$\La_Y^r(i)$ for $i=0$ and $1$ and all $r\ge r_1$, $\mu(A)=0=\mu(C)$.  Thus
$$SF(D,M^{r_1})=\mu(B)=\mu(\wtilde{L}_X,\wtilde{L}_Y).$$\qed

The transversality assumptions in Theorem \ref{yoshida} can in some contexts be
relaxed by requiring only assumptions similar to Hypothesis \ref{const}. Some care
must be taken with the steps calculating $\mu(L_3,M_3)$ and $\mu(L_9,M_9)$.

Next we give a generalization of the two theorems above by assuming the existence
of a continuously varying spectral gap.

\begin{definition}  A continuous  function $\lambda:[0,1]\to [0,\infty)$ is a {\it
spectral gap  for the family of tangential operators $S(t)$} if for each
$t\in[0,1]$,
$\lambda(t)$ is not in the spectrum of $S(t)$. 
\end{definition}

\begin{hypothesis}\label{gap} The path of tangential operators $S(t)$ has a
spectral gap
$\lambda(t)$.\end{hypothesis}

As we have remarked above, if Hypothesis \ref{gap} holds, then the decomposition
of Equation
\ref{decomp2} varies continuously.    Notice that by subdividing the path as
necessary, Hypothesis \ref{gap} can always be arranged to hold.  However, this
hypothesis by itself is not  usually sufficient to simplify the formula of Theorem
\ref{splitthm}. The following theorem gives one possible  clean  statement which
generalizes both Theorems \ref{bunke} and
\ref{yoshida}.   

Assume that Hypothesis \ref{gap} hold.  Let $A_X(t)$ and $A_Y(t)$ be continuously
varying Lagrangian subspaces of $H_{\lambda(t)}$.  Then we can take the
self-adjoint boundary conditions to be
\begin{equation}\label{apsby} B_X(t)=A_X(t)\oplus P^+_{\lambda(t)}(t)\ \
\hbox{and}  \ \  B_Y(t)=P^-_{\lambda(t)}(t)\oplus A_Y(t).\end{equation}

Then Theorem \ref{splitthm} says 
$$SF(D)=SF(D_{|X},A_X\oplus P^+_{\lambda}) + SF(D_{|Y},P^-_{\lambda}\oplus A_Y)
+\mu(A_Y(1-t),A_X(1-t)) +\sum_{i\ne3,6,9}\mu(L_i,M_i).$$ By adding hypotheses we
can make many of the extra terms vanish.

\begin{theorem}\label{3term} Assume that Hypotheses \ref{L2ker}  and \ref{gap}
hold, with spectral gap
$\lambda(t)$.   Assume  that the limiting values of extended $L^2$ solutions
$\wtilde{L}_X(i)$ and 
$\wtilde{L}_Y(i)$  are transverse for $i=0$ and $1$.   Let $A_X(t)$ and $A_Y(t)$ be
continuously varying Lagrangian subspaces of $H_{\lambda(t)}$, with
$A_Y(0)=(P^-(0)\cap H_{\la(0)})\oplus\wtilde{L}_X(0)$ and
$A_X(1)=\wtilde{L}_Y(1)\oplus(P^+(1)\cap H_{\la(1)})$. Assume further that 
$A_X(0)$ is transverse to $A_Y(0)$ and that
$A_X(1)$ is transverse to $A_Y(1)$.

Then there exists an
$r_0\ge 0$ so that for all $r\ge r_0$,
\begin{equation}\label{dobeedo} SF(D,M^r)= SF(D_{|X},A_X\oplus P^+_{\lambda}) +
SF(D_{|Y},P^-_{\lambda}\oplus A_Y) +\mu(A_X,A_Y).
\end{equation}

\end{theorem}
\proof  Since Hypothesis
\ref{L2ker} hold, 
$$\lim_{r\to\infty}\La^r_X(i)=\wtilde{L}_X(i)\oplus P^-(i)$$ and
$$\lim_{r\to\infty}\La^r_Y(i)=\wtilde{L}_Y(i)\oplus P^+(i)$$ for $i=0$ and $1$.
Since we assumed that $\wtilde{L}_X(i)$ is transverse to
$\wtilde{L}_Y(i)$ for $i=0$ and $1$, Hypothesis \ref{trans} holds, so that by
Proposition \ref{sttt}      there exists an
$r_1$ so that after replacing $M$ by $M^r$ for $r \ge r_1$, $\mu(L_1,M_1)$ and
$\mu(L_{11},M_{11})$ vanish.  

Take elliptic boundary conditions 
$B_X(t)=A_X(t)\oplus P^+_{\la(t)}(t)$  and $B_Y(t)=  P^-_{\la(t)}(t)\oplus
A_Y(t)$.  Since
   $B_Y(0)= \wtilde{L}_X(0)\oplus P^-(0)$ and $B_X(1)= 
\wtilde{L}_Y(1)\oplus P^+(1)$, Hypothesis \ref{kill} holds, so that 
  $\mu(L_i,M_i)=0$ for $i=2,5,7$, and $10$. 

The path $L_4$ is the constant path at $B_Y(1)=P^-_{\la(1)}(1)\oplus A_Y(1)$ and
$M_4$ is obtained by stretching
$\La_Y(1)$ to its adiabatic limit  $\wtilde{L}_Y(1)\oplus P^+(1)$.  Since $A_Y(1)$
is transverse to $ A_X(1)=\wtilde{L}_Y(1)\oplus (P^+(1)\cap H_{\la(1)})$ by
hypothesis,
$\mu(L_4,M_4)$ vanishes, after perhaps replacing $M$ by $M^r$ for   large enough
$r$. 

The path $L_8$ is the reverse of stretching $\La_X(0)$ to its adiabatic limit
$P^-(0)\oplus\wtilde{L}_X(0)$ and $M_8$ is the constant path at $B_X(0)=
A_X(0)\oplus P^+(0)$.  By the same argument as in the preceding paragraph  
$\mu(L_8,M_8)$ vanishes after perhaps replacing $M$ by $M^r$ for large enough
$r$. 
 
Now $\mu(B_Y(1-t),B_X(1-t))=\mu(A_Y(1-t),A_X(1-t))$.  Since $A_X(i)$ is transverse
to
$A_Y(i)$ for $i=0,1$ by hypothesis, 
$$\mu(A_Y(1-t),A_X(1-t))= \mu(A_X(t),A_Y(t)).$$

Combining  these computations  proves the theorem.\qed

The following useful corollary is just the special case of the previous theorem
when the path of tangential operators has a spectral gap $\la(t)=\epsilon$ for
$\epsilon$ small.

\begin{corollary}  Assume that Hypotheses  \ref{trans} and \ref{L2ker} hold,  that
the path of tangential operators has constant dimensional kernel, and that
$A_X(t)$ and $A_Y(t)$ are   paths of Lagrangians in
$\ker S(t)$   with 
$\wtilde{L}_X(i)=A_Y(i)$ and  $\wtilde{L}_Y(i)= A_X(i)$ for $i=0,1$.   Then  for
$r$ large enough
$$ SF(D,M^r)=SF(D_{|X},A_X\oplus P^+) + SF(D_{|Y},P^-\oplus A_Y) +\mu(A_X,A_Y).
$$
\end{corollary}
\proof Hypotheses \ref{trans} and \ref{L2ker} together imply that
$\wtilde{L}_X(i)$ is transverse to $\wtilde{L}_Y(i)$  for $i=0,1$. Thus the
hypotheses of Theorem \ref{3term} hold with $\la(t)=\epsilon$, where $\epsilon$
is smaller than the smallest non-zero eigenvalue of $S(t)$ for $t\in [0,1]$.  The
Corollary follows.\qed

We finish this subsection with a few comments about comparing Theorem
\ref{splitthm} to  Theorem C of \cite{CLM2}.   This theorem expresses the spectral
flow   as a sum of three terms; formally the theorem  looks
identical to the Formula \ref{dobeedo}, but no transversality hypotheses are 
assumed in their theorem, (although they do assume that preliminary stretching has
been done and they restrict the boundary conditions at the endpoints.)   This
might suggest that some of our $\mu(L_i,M_i)$   (in particular
$\mu(L_1,M_1)$ and $\mu(L_{11},M_{11})$) vanish  {\it without} any of the
transversality conditions.  But this is not true (examples can be concocted). The
reason that their formula  has only three terms is that their definition of
spectral flow   differs from ours in the case when transversality hypotheses do
not hold at the endpoint.   In particular, in Theorem C of \cite{CLM2}
the ``exponentially small'' eigenvalues at the endpoints of the path are treated
as if they were zero. 

To derive the  result of \cite{CLM2} from Theorem \ref{splitthm}  would require 
a more  careful analysis of the rate at which
$\La_X^r$ converges to its adiabatic limit.  An examination of Nicolaescu's proof
shows  that this rate is exponential. We speculate that by replacing the
definition of the Maslov index with the ``$1/r^2$- Maslov index'' one could
derive Theorem C of \cite{CLM2} from ours.  The article
\cite{CLM1} should be helpful for such a project.   We will not pursue this any
further since we know of no  uses for such an identification.

\subsection{Spectral flow around loops.}

One nice application of Theorem \ref{splitthm} is perhaps of more interest to index
theorists than geometric topologists.       

\begin{theorem}\label{loops} Let $D(t)$ be a loop of cylindrical, neck-compatible
Dirac operators on a manifold
$M=X\cup_\Si Y$, and let $B_X$, $B_Y$ be loops of self-adjoint elliptic boundary
conditions for the restrictions of $D$ to $X$ and $Y$ respectively.  Then 
\begin{equation}\label{RightHandSide} SF(D_{|X};B_X)+SF(D_{|Y};B_Y) +\mu(B_X,B_Y)
= 0.
\end{equation}
\end{theorem}
\proof This follows from the formula in Theorem \ref{splitthm} after much
cancellation.  First of all, the collection of all Dirac operators (on a fixed
Clifford bundle) is a vector space, hence contractible.  It follows that the
spectral flow of a loop of Dirac operators on a closed manifold is $0$.  This is
the $0$ on the right hand side of Equation \ref{RightHandSide}. 

Next, one can compute that $\mu(B_Y(1-t),B_X(1-t))=\mu(B_X(t),B_Y(t))$ if $B_X$
and $B_Y$ are loops.

It remains to show that the sum of all the other terms in Theorem
\ref{splitthm} vanish.   This is easy: the composite paths 
$$Q_1=L_1*L_2*L_4*L_5*L_7*L_8*L_{10}*L_{11}$$ and 
$$Q_2=M_1*M_2*M_4*M_5*M_7*M_8*M_{10}*M_{11}$$ are defined {\it since the path is a
loop}.   But it is immediate from the definitions of these paths that $Q_1$ is
homotopic to the constant path at $\La_X(0)$ and $Q_2$ is homotopic to the
constant path at $\La_Y(0)$.  Thus 
$$\sum_{i\ne 3,6,9}\mu(L_i,M_i)=\mu(Q_1,Q_2)=0.$$\qed

Notice that there are no hypotheses on stretching, boundary conditions, etc. in
Theorem \ref{loops}.

\subsection{Applying the method to the spectral flow on manifolds with boundary}

We conclude the user's guide  with a discussion on how to apply our method to
compute the spectral flow of the path of operators on a manifold with boundary
obtained by fixing the underlying Dirac operator but varying the boundary
conditions. 

For simplicity we consider just the special case when the boundary conditions are
of the special form $B_X(t)=A(t)\oplus P^+$ for $A_X(t)\subset \ker S$ a path of
Lagrangians;
  more general situations can be  handled by a similar  method.   This theorem
 is very similar to Theorem D of \cite{CLM2}.

\begin{theorem}\label{wbound} Let $D$ be a cylindrical, neck-compatible Dirac
operator on a smooth manifold
$X$ with boundary
$\Si$.  Let $A(t)\subset \ker S$ be a path
of Lagrangian subspaces of the kernel of the tangential operator and let
$B(t)=A(t)\oplus P^+$ be the corresponding path of elliptic boundary conditions. 
Let $\nu$ be in the nonresonance range and let
$M(t)$ be the path starting at
$\La_X$ and stretching to the adiabatic limit $P^-_\nu\oplus L_X$ (given in Lemma
\ref{stretchpath}).  Let $\wtilde{L}_X\subset \ker S$ denote the limiting values of
extended $L^2$ solutions.

Then
$$SF(D, B)= \mu(\wtilde{L}_X, A(t)) + \mu(M(t),A(0)\oplus P^+) + \mu(M(1-t),
A(1)\oplus P^+).$$

In particular, if $P^-_\nu\oplus L_X$ is transverse to $A(0)\oplus P^+$ (resp.
transverse to
$A(1)\oplus P^+$) then after replacing $X$ by $X^r$ for $r$ sufficiently large,
$\mu(M(t),A(0)\oplus P^+)=0$ (resp.\\
$
\mu(M(1-t), A(1)\oplus P^+)=0$). Hence if both transversality conditions hold, 
$$SF(D,B)=\mu(\wtilde{L}_X,A(t)).$$
\end{theorem}
\proof  First, $SF(D,B)=\mu(\La_X,B)$ by Theorem 
\ref{NicolaescusSplittingTheoremWithBoundary}.    Apply the method as follows.
\begin{enumerate}
\item Let $L_1(t)= M(t)$ and let $M_1(t)$ be the constant path at $A(0)\oplus P^+$.
So $\mu(L_1,M_1)=\mu(M(t), A(0)\oplus P^+)$.
\item Let $L_2(t)$ be the path defined in Equation \ref{rot} and let $M_2$  be the
constant path at $A(0)\oplus P^+$.  Then $\mu(L_2,M_2)=0$ by Lemma \ref{cann}. 
\item Let $L_3$ be the constant path at $P^-\oplus\wtilde{L}_X$ and let
$M_3(t)=B(t)= A(t)\oplus P^+$. So $\mu(L_3,M_3)= \mu(\wtilde{L}_X,A(t))$.
\item Let $L_4$ be $L_2$ run backwards and $M_4$ be the constant path at
$A(1)\oplus P^+$. Then $\mu(L_4,M_4)=0$ by Lemma \ref{cann}.  
\item Let $L_5$ be $L_1$ run backwards and $M_5$ the constant path at $A(1)\oplus
P^+$.
 Then $\mu(L_5,M_5)=\mu(M(1-t), A(1)\oplus P^+)$.
\end{enumerate}

Thus $L_1*L_2*L_3*L_4*L_5$ is defined and homotopic rel endpoints to the constant
path at
$\La_X$.  Also $M_1*M_2*M_3*M_4*M_5$ is defined and homotopic rel endpoints to the
path
$B$. Applying the homotopy invariance and additivity of the Maslov index finishes
the proof.\qed

\section{Concluding remarks}  

We finish with a few comments about  Theorem \ref{splitthm}.  First, there is a
certain asymmetry in   the formula   with respect to the roles that $X$ and $Y$
play. This turns out to be useful sometimes, for example Hypothesis \ref{kill}
only restricts
$B_Y$ at one endpoint and $B_X$ at the other, rather than restricting both at each
endpoint.   

Another comment is that the sums $\mu(L_2,M_2)+\mu(L_7,M_7)$ and
$\mu(L_5,M_5)+\mu(L_{10},M_{10})$ (the terms depending on  the
auxiliary choice of the paths
$L_2$ and
$M_5$) depend only on the endpoints of these paths.  Thus each of these sums
could be thought of as  a single quantity, and perhaps expressed in terms of
invariants (such as the Maslov triple index) of the endpoints alone, without
making any reference to the choice of $L_2$ and $M_5$.  As we have seen, it is
nevertheless convenient for calculation to have the formula expressed the way we
did.

Last (but not least), one significant benefit of our formulation is that since
our formula expresses the spectral flow entirely as a sum of Maslov indices, with
the {\it ordered} pairs $(L_i,M_i)$ explicitly described, it is easy to keep the
signs and conventions under control when carrying out spectral flow
calculations.

 \vskip4ex

\begin{appendix}
\section{The proof of Lemma
\ref{stretchpath}, by K. P. Wojciechowski}

The set up is as follows. We are given a Dirac operator $D$ on  manifold $X$
with boundary    in cylindrical form $D=J(\partial/ \partial u + S)$ on a collar
$\Si\times[-1,0]$ of the boundary $\Si=\Si\times \{0\}$. This extends to an
operator on $X^r=X\cup \Si\times [0,r]$ in the obvious way.   To this extension
we associate the Cauchy data spaces
$\La^r$.

Nicolaescu's adiabatic limits theorem, Theorem \ref{NicoAdiabaticLimitTheorem},
says that the path (with $r(t)=1/(1-t)$)
$$t\mapsto \begin{cases}
\La^{r(t)}& t<1,\\ P^-_{\nu_0}\oplus  L_X(D) &t=1.
\end{cases}
$$  is continuous at $t=1$.  What must be shown is that this path  is continuous
at finite neck lengths $r$, that is, that the Cauchy data spaces $\La^r$ vary
continuously in $r$.  Continuity is measured in the gap topology, or
equivalently in the norm of the associated projections.

For notational convenience we will prove continuity at $r=0$; by
reparameterizing continuity at all $r$ follows easily.
\vskip4ex

Let $\nu$ be a number in the nonresonance range for $D$ on $X=X^0$.
Thus $\La^0\cap P^+_{\nu'}=0$ for all $\nu'\ge \nu$.   We will make frequent use
of the splitting $L^2(E_{|\Si})= P^-_\nu\oplus H_\nu\oplus P^+_\nu$. Notice that
the tangential operator $S$ preserves this splitting since the summands are
defined by the eigenspace decomposition of $S$.  The almost complex structure
$J$ of Equation \ref{complx} preserves $H_\nu$ and interchanges $P^+_\nu$ and
$P^-_\nu$.

 We will often use   the fact that if $\alpha$ is a section of $E$ on
the cylinder $\Si\times [-1,r]$ which satisfies $D\alpha=0$, then writing
$\alpha_{|\Si\times\{u\}}=\alpha(u)$ for $u\in [-1,r]$,
$$\alpha(u)=e^{(t-u)S}\alpha(t).$$

Let $L=\proj_{H_\nu}(\La^0\cap(H_\nu\oplus P^+_\nu))$. Clearly $JL\oplus
P^+_\nu$ is transverse to $\La^0$.

\begin{lemma}\label{tralem} For each $r\ge -1$, $\La^r$ is transverse to
$(e^{-rS}JL)\oplus P^+_\nu$.
\end{lemma}
\proof   Suppose that $v\in \La^r\cap( (e^{-rs}JL)\oplus P^+_\nu)$. Then there
exists $\al $ a section of $E$ on $X^r$ so that $D\al=0$ and the
restriction of $\al$ to
$\Si\times\{r\}$ equals $0$. Thus $\al(u)=e^{(r-u)S}v$ for $u\in[-1,r]$.  But
this formula defines an extension of $\al(u)$ for all $u\in[-1,\infty)$ since
$H_\nu$ is finite dimensional and since the restriction of $e^{(r-u)S}$ to
$P^+_\nu$ exponentially decays as $u\to \infty$.  Hence $\al$ extends to a
bounded smooth section on $X^u$ for all $u>-1$, and the extension  satisfies
$D\al=0$. In particular, $\al(0)$ is defined, equals $e^{rS}v$, and lies in
$\La^0$. Since $v\in (e^{-rS}JL)\oplus P^+_\nu$, $\al(0)\in e^{(r-0)S}
((e^{-rS}JL)\oplus P^+_\nu)= JL\oplus P^+_\nu$. By the choice of $L$ this
implies that $\al(0)=0$ and so also $v=0$.\qed

For convenience we introduce some notation for certain projections
$L^2(E_{|\Si})$.

\begin{enumerate}
\item The {\it orthogonal Calderon projection} $\cP^r:L^2(E_{|\Si})\to
L^2(E_{|\Si})$ is the orthogonal projection to the Cauchy data space $\La^r$.

\item The {\it negative spectral projection} $\pi_-:L^2(E_{|\Si})\to
L^2(E_{|\Si})$ is the orthogonal projection to the space $P^-(S)$, the negative
eigenspan of the tangential operator $S$.

\item Fix $\nu\ge 0$ and suppose $L\subset H_\nu$ is a Lagrangian (thus
$P^-_\nu\oplus L$ is a Lagrangian in $L^2(E_{|\Si}))$.  Define $\pi_{-,L}
:L^2(E_{|\Si})\to  L^2(E_{|\Si})$ to be the orthogonal projection to
$P^-_\nu\oplus L$.
\end{enumerate}

What must be shown is that the projections $\cP^r$ are continuous in norm as $r$
varies.  It follows from the results in Chapter 12 and 14 of \cite{BW} that
$\cP^r$, $\pi_-$, and
$\pi_{-,L}$ are pseudodifferential of order $0$.

Let $L$ and $\nu$ be as in Lemma \ref{tralem}. For notational ease, define
$$M_r=e^{rS}L\subset H_\nu$$ and $$\pi_r=\pi_{-,M_r}:L^2(E_{|\Si})\to M_r\oplus
P^-_\nu.$$ Notice that $M_r$ varies continuously in $r$, and hence so does
$\pi_r$.  The difference $\pi_- -\pi_r$ has image in $H_\nu$, a finite
dimensional space of smooth sections, and hence is a smoothing operator.
Corollary  14.3 of
\cite{BW} shows that the pseudo-differential operators $\cP^r$ and $\pi_-$ have
the same principal symbol. Putting these facts together shows that $\cP^r-\pi_r$
is a pseudo-differential operator of order at most $-1$, and in particular is a
compact operator.

\begin{lemma}\label{project} The restriction of $\pi_r$ to $\La^r$ induces an
isomorphism
$$\pi_r:\La^r\to M_r\oplus P^-_\nu$$
\end{lemma}

\proof It is easy to observe that the operator
$$\pi_r:\La^r\to M_r\oplus P^-_\nu$$
is a Fredholm operator (see \cite{BW}), hence in particular it has closed
range.
The kernel of this map is $\La^r\cap (M_r\oplus P^-_\nu)^\perp=
\La^r\cap (JM_r\oplus P^+_\nu)=0$. 
Since $\La^r$ and $M_r\oplus P^-_\nu$ are Lagrangians  the isometry $J$
identifies the  cokernel with the kernel of $\pi_r$, and so the map is
surjective.\qed

Since the map of Lemma \ref{project} is an isomorphism,    the Cauchy data space
$\La^r$ can be expressed as  a graph of a bounded operator
$$k_r: M_r\oplus P^-_\nu\to JM_r\oplus P^+_\nu;$$
here $k_r$ is the composite of the inverse of $\pi_r:\La^r\to M_r\oplus
P^-_\nu$ and the orthogonal projection to $JM_r\oplus P^+_\nu$.

Hence
\begin{equation}
\La^r=\{ (v, k_r(v))\ | \ v\in M_r\oplus P^-_\nu\}.
\end{equation}

\vskip3ex

Let $r>-1$. Choose   $v_-\in M_r\oplus P^-_\nu$. Hence $v=v_- + k_r(v_-)\in
\La^r$.   Thus there exists a section  $\alpha$ in $\ker D$ on $X^r$  with
$\al(r)=v$.  As observed above,  on the cylinder $\Si\times[-1,r]$ $\alpha$ has
the form
\begin{equation}\label{expeq1}
\alpha(u)= e^{(r-u)S}v=e^{(r-u)S}v_- + e^{(r-u)S}k_r(v_). \end{equation}

On the other hand,  for $u\leq r$  $\alpha(u)\in \La^u$,  and  taking $u=-1$ we
have
\begin{equation}\label{expeq2}
\alpha(-1)=w_-+ k_{-1}(w_-)\end{equation} for some $w_-\in M_{-1}\oplus
P^+_\nu$.

Combining Equations \ref{expeq1} and \ref{expeq2} yields
$$w_-= e^{(r+1)S}v_- \ \ \hbox{ and } e^{-(r+1)S}k_{-1}(w_-)=k_r(v_-).$$ hence
\begin{equation}\label{expeq3} k_r=e^{-(r+1)S_+}k_{-1}e^{(r+1)S_-}
\end{equation} where we have denoted the restriction of $S$ to $H_\nu\oplus
P^\pm_\nu$ by $S_\pm $ for clarity.

For the next lemma, we recall the standard fact  that the
operators $e^{tS_-}:H_\nu
\oplus P^-_\nu\to H_\nu
\oplus P^-_\nu$  and
$e^{-tS_+}:H_\nu\oplus P^+_\nu\to H_\nu\oplus P^+_\nu $ are norm-continuous in
$t$ for $t$ away from $0$.  

\begin{lemma}\label{estimate}
$$\lim_{r\to 0}\norm k_r-k_0\norm =0.$$
\end{lemma}
\proof Using Equation \ref{expeq3} we compute, for $-\frac{1}{2} \leq r
\leq\frac{1}{2}$
\begin{eqnarray*}
\norm k_r -k_0\norm&=&
\norm  e^{-(r+1)S_+}k_{-1}e^{(r+1)S_-}- e^{-S_+}k_{-1}e^{S_-}\norm  \\
&\leq&
\norm (e^{-(r+1)S_+}k_{-1}e^{(r+1)S_-}-
     e^{-(r+1)S_+}k_{-1}e^{(r+1)S_-}\norm  \\
  &&\hskip1in +\ \norm (e^{-(r+1)S_+}k_{-1}e^{S_-}-
     e^{-S_+}k_{-1}e^{S_-}\norm \\ &\leq& \norm e^{-(r+1)S_+}\norm\ \norm
k_{-1}(e^{(r+1)S_-}- e^{S_-})\norm +
 \norm (e^{-(r+1)S_+}-e^{-S_+}) k_{-1}\norm\ \norm e^{S_-}\norm \\
&\leq&
C_1
\norm k_{-1}(e^{(r+1)S_-}- e^{S_-})\norm + C_2  \norm (e^{-(r+1)S_+}-e^{-S_+})
k_{-1}\norm.
\end{eqnarray*} The last inequality follows from the continuity of
$e^{(r+1)S_+}$ in norm for $r$ near $0$ and the fact that $e^{S_-}$ is
independent of $r$.  Continuing the estimate using the fact that $k_{-1}$ is
bounded, we obtain
\begin{equation}
\norm k_r -k_0\norm \leq  \norm k_{-1}\norm
\left( C_1 \norm e^{(r+1)S_-}-e^{S_-}\norm + C_2 \norm e^{-(r+1)S_+} -
e^{-S_+}\norm\right).
\end{equation}
The right hand side approaches $0$ as $r\to 0$ since $e^{tS_-}$ and $e^{-tS_+}$
are continuous in norm at $t=1$. This proves the Lemma.\qed

In the decomposition
$$L^2(E_{|\Si})= (M_r\oplus P^-_\nu)\oplus(JM_r\oplus
P^+_\nu)$$ the matrix
$$Q_r=\begin{pmatrix} 1&0\\
                     k_r&0\end{pmatrix} $$
is  a (non-orthogonal) projection to $\La^r$.  The formula of Lemma 12.8 in
\cite{BW}  shows that

 \begin{eqnarray}\label{explic}
\cP^r&=&Q_rQ_r^*\left( Q_rQ_r^* + (Id - Q^*_r)(Id-Q_r)\right)^{-1} \nonumber\\
&=&\begin{pmatrix}
(Id + k_r^*k_r)^{-1} & (Id + k_r^*k_r)^{-1}k_r^* \\ k_r(Id + k_r^*k_r)^{-1}
& k_r(Id
+ k_r^*k_r)^{-1}k_r^*
\end{pmatrix} 
\end{eqnarray}

Thus, the $\cP^r$ are continuous in $r$ for $r$ near $0$,
completing the proof of Lemma
\ref{stretchpath}.\qed
\vskip4ex

We note that the proof of Lemma \ref{estimate} can can be  modified to show
that
 $\lim_{r\to
\infty}\La^r=(\lim_{r\to \infty}M_r)\oplus P^-_\nu=(\lim_{r\to
\infty}e^{rS}L)\oplus P^-_\nu$. In contrast to the proof of continuity at
finite $r$ given above, in this case one  must be careful with the estimates
over the finite-dimensional piece $H_\nu$.  The argument is straightforward,
and  is essentially the proof given in  Nicolaescu
\cite{Nico}.

Let us also notice that in fact we proved here the following important result.

\begin{theorem}\label{smooth}
The difference  $\cP^r - \pi_-$ is an operator with a
smooth kernel
\end{theorem}

\begin{proof}
The difference $\pi_r - \pi_-$ is a smoothing operator and using
Equation
\ref{explic}
the difference $\cP^r - \pi_r$ can
be represented as
$$\cP^r - \pi_r = \begin{pmatrix}
(Id + k_r^*k_r)^{-1} - Id & (Id + k_r^*k_r)^{-1}k_r^* \\ k_r(Id +
k_r^*k_r)^{-1} & k_r(Id
+ k_r^*k_r)^{-1}k_r^*
\end{pmatrix} .$$

All entries in the formula presented above are smoothing operators due to
the fact that $k_r$ has a
smooth kernel.
\end{proof}

S. Scott proved this Theorem in the non-resonant case (see
\cite{Sc}). The proof given above  basically extends his proof to cover the
general case. A different proof, purely analytical, was offered by G. Grubb in
(see \cite{Gr}).

\end{appendix}

\parskip 1.75\parskip plus 3pt minus 1pt 
\renewcommand{\baselinestretch}{1}

\vskip4ex
{\obeylines  {\em
Mark Daniel
Advanced Power Technologies, Inc.
Washington, DC 20037
e-mail: amdaniel@apti.com
\vskip3ex

Paul Kirk
Indiana University
Bloomington IN, 47405
e-mail: pkirk@indiana.edu

\vskip3ex
Krzysztof P. Wojciechowski
Indiana University-Purdue University at Indianapolis
Indianapolis, IN 46202
e-mail: kwojciechowski@math.iupui.edu
}}
\vfill\eject

\bibliographystyle{amsplain}

\end{document}